\documentclass[11pt]{amsart}

\usepackage{graphicx,amssymb}
\theoremstyle{definition}
\newtheorem{de}{Definition}[section]
\theoremstyle{plain}
\newtheorem{lm}[de]{Lemma}
\newtheorem{pr}[de]{Proposition}
\newtheorem{te}[de]{Theorem}

\theoremstyle{remark}
\newtheorem{re}[de]{Remark}

\newtheorem{ex}[de]{Example}
\def\bea{\begin{eqnarray*}}
\def\eea{\end{eqnarray*}}

\font\nb=msbm10

\def \Z{\hbox{\nb Z}}

\def \Q{\hbox{\nb Q}}

\def\id{\textrm{Id}}

\def\ot{\otimes}

\def\phi{\varphi}
\def\ot{\otimes}

\newcommand{\nbpt}{18}
\newcommand{\figtotext}[3]{\begin{array}{c}\includegraphics[width=#1pt,height=#2pt]{#3}\end{array}}
\newcommand{\overcrossing}{\figtotext{\nbpt}{\nbpt}{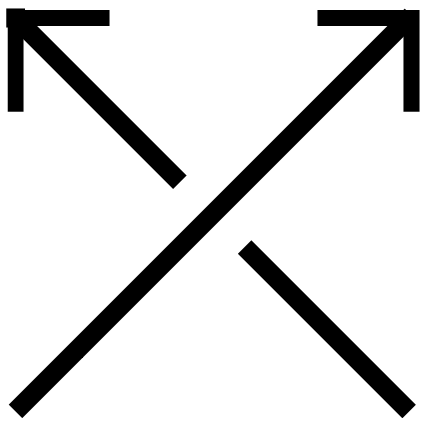}}
\newcommand{\undercrossing}{\figtotext{\nbpt}{\nbpt}{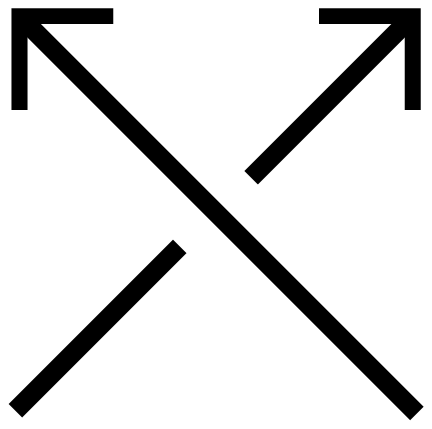}}
\newcommand{\smoothing}{\figtotext{\nbpt}{\nbpt}{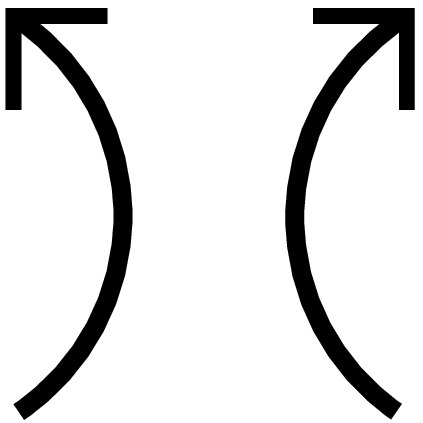}}
\newcommand{\unknot}{\figtotext{\nbpt}{\nbpt}{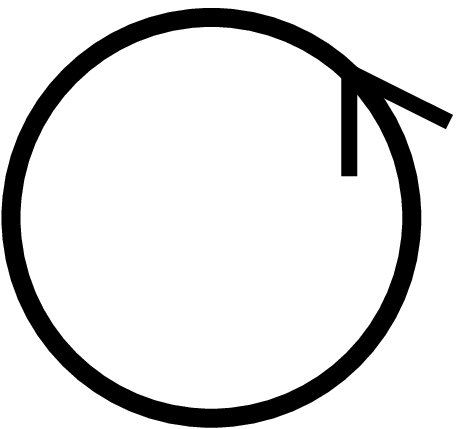}}

\begin{document}
\title[]{Yang--Baxter operators arising from algebra structures
and the Alexander polynomial of knots}
\author[G. Massuyeau]{Gw\'ena\"el Massuyeau}
\author[F. F. Nichita]{Florin F. Nichita}
\date{January, 2004}
\begin{abstract}
In this paper, we consider the problem of constructing knot invariants
from Yang--Baxter operators associated to algebra structures. We first
compute the enhancements of these operators. Then, we conclude that Turaev's
procedure to derive knot invariants from these enhanced operators, as
modified by Murakami, invariably produces the Alexander polynomial of
knots.

\end{abstract}
\maketitle

\section{Introduction}
The Yang--Baxter equation and its solutions, the Yang--Baxter operators,
first appeared in theoretical physics and statistical mechanics.
Later, this equation has emerged in other fields of mathematics such 
as quantum group theory.
Some references on this topic are \cite{Kassel,LR}.

The Yang--Baxter equation also plays an important role in knot theory.
Indeed, Turaev has described in \cite{T}
a general scheme to derive an invariant of oriented links
from a Yang--Baxter operator, provided this one can be ``enhanced''.
The Jones polynomial \cite{Jones} and its two--variable extensions,
namely the Homflypt polynomial \cite{HOMFLY,PT} and the Kauffman 
polynomial \cite{K},
can be obtained in that way by ``enhancing'' some Yang--Baxter operators
obtained in \cite{Jimbo}. Those solutions of the Yang--Baxter 
equation are associated to simple Lie algebras
and their fundamental representations. The Alexander polynomial can be derived
from a Yang--Baxter operator as well, using a slight modification of 
Turaev's construction \cite{M}.

More recently, D\u asc\u alescu and Nichita have shown in \cite{DN} 
how to associate a Yang--Baxter operator
to any algebra structure over a vector space, using the associativity 
of the multiplication.
This method to produce solutions to the Yang--Baxter equation, 
initiated in \cite{N}, is quite simple.

In this paper, we consider the problem of applying Turaev's method to
Yang--Baxter operators derived from algebra structures. In general,
finding the enhancements of a given Yang--Baxter operator can be 
difficult or lengthy.
In the case of Yang--Baxter operators associated to algebra structures,
the simplicity of their definition makes the search for enhancements 
an easy task.
We do this here in full generality.
We conclude from this computation that the only invariant which can 
be obtained from
those Yang--Baxter operators is the Alexander polynomial of knots.
Thus, in a way, the Alexander polynomial is the knot invariant corresponding to
the axioms of (unitary associative) algebras.
Note that specializations of the Homflypt polynomial had to be expected
from those Yang--Baxter operators since they have degree $2$  minimal 
polynomials.\\

The paper is organized as follows. In \S \ref{sec:YB},
we recall how to associate to any (unitary associative) algebra
a Yang--Baxter operator. Next, in \S \ref{sec:invariants},
we review  Turaev's procedure to derive a knot invariant
from a Yang--Baxter operator as soon as this one can be enhanced.
Following Murakami \cite{M}, we recall how this method can be improved
in the case when the Yang--Baxter operator satisfies a certain 
redundancy property.
The fact that their minimal polynomials are quadratic implies that
the Yang--Baxter operators associated to algebra structures are redundant
(Proposition \ref{prop:quadratic}). In \S 
\ref{sec:invariants_associated_to_algebras},
we compute the enhancements of the Yang--Baxter operator associated 
to a given algebra structure.
Finally, we conclude from this calculation that Turaev's procedure, 
as modified by Murakami,
invariably produces from any of those enhancements the Alexander 
polynomial of knots
(Theorem \ref{th:conclusion}).\\

Throughout the paper, we shall use the following conventions and notations.
The letter $k$ will denote a fixed field with characteristic char$(k)\neq 2$.
For a $k$--vector space $V$ and an integer $n\geq 1$, we denote by 
$V^{\otimes n}$ the $n$--times
tensor product $V\otimes \cdots \otimes V$ over $k$. The identity map 
$V \to V$ will be denoted
by $\id_V$, or simply by $\id$ when the space $V$ is clear from the 
context. If $V$ is finite--dimensional with basis
$e=(e_1,\dots,e_d)$ and if $f: V^{\ot n}\to V^{\ot n}$ is a $k$--linear map,
we denote by $f^{j_1 \cdots  j_n}_{i_1\cdots i_n} $
the matrix element of $f$ with respect to the basis $e^{\ot n}$:
$$ f\left(e_{i_1} \ot \cdots \ot e_{i_n}\right)
= \sum_{1\leq j_1, \dots , j_n \leq d} f^{j_1 \cdots j_n}_{i_1\cdots 
i_n} \cdot e_{j_1}  \ot \cdots \ot e_{j_n}.$$
\vspace{0.1cm}

\section{Yang--Baxter operators arising from algebra structures}

\label{sec:YB}

In this section, we recall results from \cite{DN} and we follow the 
terminology used there.
For a $k$--vector space $V$, the \emph{flip} map $T:V \ot  V 
\rightarrow V \ot  V$
is defined by $T(v \ot w) = w \ot v.$ Given a $k$--linear map
$R: V \ot V \rightarrow V \ot V$,  we will write $ {R^{12}}= R \ot 
\id, \ {R^{23}}= \id \ot R$ and
${R^{13}}=(\id\ot T)(R\ot \id)(\id\ot T).$

\begin{de}
An invertible  $ k$--linear map  $ R : V \ot V \rightarrow V \ot V $
is called a \emph{Yang--Baxter operator} (or simply a \emph{YB 
operator}) if it satisfies the  equation
\begin{equation}
\label{ybeq}
R^{12}  \circ  R^{23}  \circ  R^{12} = R^{23}  \circ  R^{12}  \circ  R^{23}.
\end{equation}
\end{de}

\begin{re}
Equation (\ref{ybeq}) is usually called the \emph{braid equation}.
An operator $R$ satisfies (\ref{ybeq}) if and only if
$R\circ T$ satisfies the \emph{quantum Yang--Baxter equation}
$$
R^{12}  \circ  R^{13}  \circ  R^{23} = R^{23}  \circ  R^{13}  \circ  R^{12}.
$$
\end{re}

Given a (unitary associative) $k$--algebra $A$ and scalars $x,y,z 
\in k$, we consider
the $k$--linear map $R_{x,y,z}: A \ot A \rightarrow A \ot A $ defined by
$$R_{x,y,z}(a \ot b)= x\cdot ab \ot 1 + y\cdot 1 \ot ab - z\cdot a \ot b$$
for any $ a,b \in A $. The following theorem determines
the situations where this map is a Yang--Baxter operator.

\begin{te}[D\u asc\u alescu--Nichita \cite{DN}]
\label{th:DN}
Let A be a $k$--algebra of dimension at least $2$ and let $ x,y,z 
\in k $ be scalars. Then
$R_{x,y,z}$ is a Yang--Baxter operator if and only if one of the 
following conditions holds:
\begin{itemize}
\item[(i)] $ x=z \neq 0, \ y \neq 0 $,
\item[(ii)] $ y=z \neq 0, \ x \neq 0 $,
\item[(iii)] $ x=y=0 , \ z \neq 0 $.
\end{itemize}
If so, its inverse is given by 
${R_{x,y,z}}^{-1}=R_{y^{-1},x^{-1},z^{-1}}$ in cases {\rm{(i)}} and 
{\rm{(ii)}},
and by ${R_{0,0,z}}^{-1}=R_{0,0,z^{-1}}$ in case {\rm{(iii)}}.
\end{te}

\noindent
We will be only interested in the non--trivial cases (i) and (ii) 
(which are equivalent through inversion)
and we will denote by $R_{x,y}$ the operator $R_{x,y,x}$. Observe that
\begin{equation}
\label{eq:quadratic}
R_{x,y} - xy \cdot {R_{x,y}}^{-1}=(y-x)\cdot \id^{\ot 2}
\end{equation}
and that, in particular, the minimal polynomial of $R_{x,y}$ is of degree $2$.

Next result gives a necessary and sufficient condition for two 
Yang--Baxter operators
of that kind to be isomorphic.

\begin{de}
\label{def:isoYB}
Two $k$--linear maps $ R: V \ot V \rightarrow V \ot V $ and
$ R': V' \ot V' \rightarrow V' \ot V' $ are
said to be \emph{isomorphic} if there exists an invertible $k$--linear map
$ f: V \rightarrow V' $ such that $ R' \circ (f \ot f) = (f \ot f) \circ R $.
\end{de}

\begin{pr}[D\u asc\u alescu--Nichita \cite{DN}]
Let $A$ and $A'$ be $k$--algebras of dimension at least $2$, and let 
$x,y,x',y'$ be
nonzero scalars. Then the YB operators $R_{x,y}$ and $R_{x',y'}$
associated to $A$ and $A'$ respectively are isomorphic if, and only 
if, $ x=x',\ y=y'$
and the $k$--algebras $A$ and $A'$ are isomorphic.
\end{pr}

\begin{re}
\label{rem:isoYB}
In fact, according to the proof of \cite[Proposition 3.1]{DN},
a bijective $k$--linear map $f: A \to A'$ is such that $ R_{x',y'} 
\circ (f \ot f) = (f \ot f) \circ R_{x,y}$
if and only if $ x=x', y=y'$ and $f$ is, up to multiplication by a scalar,
an isomorphism of $k$--algebras. (Note that the assumption 
char$(k)\neq 2$ is used here.)
\end{re}

\section{Invariants of oriented links derived from Yang--Baxter operators}

\label{sec:invariants}

In this section, we shortly review Turaev's procedure \cite{T} to 
derive an invariant of oriented links
from a Yang--Baxter operator, provided one can ``enhance'' it.
We also recall a slight modification of this construction, due to 
Murakami \cite{M}, which applies
to Yang--Baxter operators verifying a certain redundancy property.

\subsection{Enhanced Yang--Baxter operators}

\label{subsec:enhanced_YB}

Let $V$ be a finite--dimensional $k$--vector space
and let $f:V^{ \ot {n}}\to V^{ \ot {n}}$ be a $k$--linear map. We 
pick a basis $e=(e_1,\dots, e_d)$ of $V$ and,
following our convention, we denote $f^{j_1 \cdots  j_{n-1} 
j_n}_{i_1\cdots i_{n-1} i_n} $
the matrix element of $f$ with respect to the basis $e^{\ot n}$.
Then, the \emph{operator trace} of $f$ is the $k$--linear map
$\hbox{Sp}_n(f):V^{ \ot {n-1}}\to V^{ \ot {n-1}}$ defined by
$$ \hbox{Sp}_n(f)\left(e_{i_1} \ot \cdots \ot e_{i_{n-1}}\right)
= \sum_{\substack 1\leq j_1, ... , j_{n-1}, l \leq d}
f^{j_1\cdots j_{n-1} l}_{i_1\cdots i_{n-1}  l} \cdot e_{j_1} \ot 
\cdots \ot e_{j_{n-1}}.$$

\begin{de}
\label{def:enhancement}
An \emph{enhanced Yang--Baxter operator} is a quadruplet
$ S = ( R, \mu, \alpha, \beta ) $ where
$ R: V\ot V\to V\ot V$ is a YB operator on a finite--dimensional 
$k$--vector space $V$,
$\mu: V\to V$ is a $k$--linear bijective map and $ \alpha, \beta$ are 
nonzero scalars in $k$ satisfying
\begin{itemize}
\item[$\left(E_1\right)$] $R \circ ( \mu \ot \mu ) = ( \mu \ot \mu ) \circ R$,
\item[$\left(E_2^\pm\right)$] $\hbox{Sp}_2\left( R^{\pm 1} \circ ( 
\id \ot \mu )\right) = \alpha^{\pm 1} \beta \cdot \id.$
\end{itemize}
\end{de}

For any integer $n\geq 1$, let $B_n$ denotes the $n$--string braid 
group which can be presented as
$$
B_n =\left\langle \sigma_1,\dots, \sigma_{n-1}\ \vert\ \sigma_i 
\sigma_{i+1} \sigma_i = \sigma_{i+1} \sigma_i \sigma_{i+1},
\sigma_i \sigma_j = \sigma_j \sigma_i \ \textrm{if}\ |i-j|>1  \right\rangle.
$$
Any Yang--Baxter operator $R$ defines a linear representation
$$
\rho_R : B_n \longrightarrow \hbox{End} \left(V^{\ot n}\right)
$$
of the braid group by putting
$$
\rho_R\left(\sigma_i\right)= \id^{\ot i-1} \ot R \ot \id^{\ot n-i-1}.
$$

If the Yang--Baxter operator $R$ can be enhanced to 
$S=(R,\mu,\alpha,\beta)$, then one can go further
and define an invariant of oriented links as follows. First,
any braid $b\in B_n$ leads to an oriented link $\widehat{b}$ by 
\emph{closing} its $n$ strings.
A theorem of Alexander asserts that any link is the closure of a 
braid, and a theorem of Markov
gives a necessary and sufficient condition for two braids (with 
possibly different numbers of strings)
to have isotopic closures. Next, we can associate to any braid $b\in 
B_n$ its \emph{Markov trace}
$$
T_S (b) = { \alpha }^{-w(b)} \beta^{-n}\cdot
\hbox{Trace}\left(\rho_R(b) \circ \mu^{\ot n} \right) \in k,
$$
where $w: B_n \to \Z$ is the group homomorphism defined by 
$w(\sigma_i)=1$ for any $i=1,\dots,n-1$.

Next theorem is proved from Markov's theorem, the definition of an
enhanced YB operator and properties of the trace.

\begin{te}[Turaev \cite{T}] Let $S$ be an enhanced YB operator.
If two braids $b_1$ and $b_2$ have isotopic closures, then their 
Markov traces $T_S(b_1)$
and $T_S(b_2)$ are equal. So, there is an isotopy invariant $X_S$ of 
oriented links defined
by $X_S\left(\widehat{b}\right)=T_S(b)$.
\end{te}

\subsection{Redundant Yang--Baxter operators}

\label{subsec:redundant_YB}

If $R\in \hbox{End}\left(V^{\ot 2}\right)$ is a Yang--Baxter 
operator, we define
$A_{R,n}$ to be the subspace of $\hbox{End}\left(V^{\ot n}\right)$ 
generated by the image $ \rho_R (B_n)$.
In particular, $A_{R,n}$ is a subalgebra of $\hbox{End}\left(V^{\ot n}\right)$
and we have $ A_{R,1} = \langle \id \rangle \subset \hbox{End}(V)$.

\begin{de}
A \emph{redundant Yang--Baxter operator} is an enhanced YB operator
$ S = ( R, \mu, \alpha, \beta ) $ such that
$$
\forall x \in A_{R,n}, \
\hbox{Sp}_n^\mu(x) \in A_{R, n-1}.
$$
Here, for any $x\in \hbox{End}\left(V^{\ot n}\right)$, 
$\hbox{Sp}_n^\mu(x)$ denotes
$\hbox{Sp}_n \left(x \circ \left(\id^{\ot n-1} \ot \mu\right)\right)$.
\end{de}

Let us observe that, under the assumption that $S$ is redundant,
the endomorphism of $V$ associated to a braid $b\in B_n$ by
$$
\hbox{Sp}_2 \cdots \hbox{Sp}_n\left(\rho_R(b) \circ\left(\id \ot 
\mu^{\ot n-1}\right)\right)
= \hbox{Sp}_2^\mu \cdots \hbox{Sp}_n^\mu\left(\rho_R(b)\right)
$$
is a multiple of the identity and, so, can be regarded as a scalar.
This allows us to define the \emph{modified Markov trace} of $b\in B_n$ to be
$$
T_{S,1}(b)= \alpha^{-w(b)} \beta^{-n}\cdot
\hbox{Sp}_2  \cdots \hbox{Sp}_n\left(\rho_R(b) \circ\left(\id \ot 
\mu^{\ot n-1}\right)\right)\in k.
$$
In this case, we will have that
\bea
T_S(b) &=& { \alpha }^{-w(b)} \beta^{-n}\cdot 
\hbox{Trace}\left(\rho_R(b) \circ \mu^{\ot n} \right)\\
&=& { \alpha }^{-w(b)} \beta^{-n}\cdot
\hbox{Sp}_1 \hbox{Sp}_2 \cdots \hbox{Sp}_n\left(\rho_R(b) \circ 
\mu^{\ot n} \right)\\
&=& { \alpha }^{-w(b)} \beta^{-n}\cdot
\hbox{Sp}_1^\mu \hbox{Sp}_2^\mu \cdots \hbox{Sp}_n^\mu\left(\rho_R(b)\right)\\
&=& \hbox{Trace}(\mu)\cdot T_{S,1}(b).
\eea
\begin{te}[Murakami \cite{M}]
Let $S$ be a redundant YB operator.
If two braids $b_1$ and $b_2$ have isotopic closures, then their 
modified Markov traces
$T_{S,1}(b_1)$ and $T_{S,1}(b_2)$ are equal.
So, there is an isotopy invariant $X_{S,1}$ of oriented links defined
by $X_{S,1}\left(\widehat{b}\right)=T_{S,1}(b)$.
\end{te}
\begin{re}
\label{rem:trivialize}
If $S=(R,\mu,\alpha,\beta)$ is a redundant Yang--Baxter operator, 
then we have that
$X_{S}\left(L\right)=\hbox{Trace}(\mu)\cdot X_{S,1}(L)$ for any 
oriented link $L$.
In particular, if Trace$(\mu)$ vanishes, then $X_S$ does too:
this is the case when it is worth working with modified Markov traces.
\end{re}

Next proposition will be proved thanks to arguments from 
\cite[Appendix A]{M}. There,
Murakami shows redundancy of a particular (colored) enhanced 
Yang--Baxter operator,
from which he derives the (multivariable) Alexander polynomial.

\begin{pr}
\label{prop:quadratic}
Let $ S = ( R, \mu, \alpha, \beta ) $ be an enhanced Yang--Baxter 
operator such that
$R^2 + b\cdot R + c \cdot{\rm{Id}}^{\ot 2} = 0 \in 
{\rm{End}}\left(V^{\ot 2}\right) $, where $c \in k^* $ and $ b \in k 
$.
Then, $S$ is redundant.
\end{pr}
\noindent
In order to prove this proposition, we will need the following lemma.
(Compare it with \cite[Lemma A.2]{M}.)

\begin{lm}
Let $R: V\ot V\to V\ot V$ be a YB operator such that
$R^2 + b\cdot R + c \cdot{\rm{Id}}^{\ot 2} = 0 \in 
{\rm{End}}\left(V^{\ot 2}\right) $, where $c \in k^* $ and $ b \in k 
$.
Then, for any integer $n\geq 2$, we have that
$$ A_{R,n}= A_{R,n-1}+ \left\langle A_{R,n-1}\cdot ({\rm{Id}}^{\ot 
n-2} \ot R) \cdot A_{R,n-1} \right\rangle$$
where $A_{R,n-1}$ is regarded as a subspace of $A_{R,n}$ under the 
natural inclusion $B_{n-1}\subset B_n$.
\end{lm}
\begin{proof}
The proof goes by induction on $ n \geq 2 $.

Suppose that $n=2$. Since $A_{R,1} = \langle \id \rangle$,  we are 
reduced to verify that
$A_{R,2} = \left\langle \id^{ \ot 2}\right\rangle + \langle R 
\rangle$. This holds true since,
for any integer $m$, $\rho_R\left({\sigma_1}^m\right)$ equals $R^m$. But, $R^m$
can be expressed as a linear combination of Id and $R$, because the 
minimal polynomial of $R$ is quadratic.

Let us now suppose that the lemma is satisfied at rank $n-1$. We want 
to prove it at rank $n$. Denote
$$
C_n= A_{R,n-1}+ \left\langle A_{R,n-1}\cdot ({\rm{Id}}^{\ot n-2} \ot 
R) \cdot A_{R,n-1} \right\rangle.
$$
We certainly have that $ \rho_R( b) \in C_n $ for any $ b \in B_{n-1} 
\subset B_n $ and for
$b = \sigma_{n-1}$. So, it is enough to prove the following\\

\noindent
\emph{Claim.}
If $ b, c \in B_n $ are such that $ \rho_R( b) \in C_n $ and
$ \rho_R(c) \in C_n $, then  we have that $ \rho_R( bc) \in C_n $ too.\\

\noindent

Let us prove that $ \rho_R( bc)= \rho_R( b) \rho_R(c)$
belongs to $C_n$ knowing that $\rho_R(b)$ and $\rho_R(c)$ do.
Since $A_{R,n-1}$ is a subalgebra of $A_{R,n}$,
we can suppose with no loss of generality that
$$
\rho_R( b) = b_1 \left(\id^{\ot n-2} \ot R\right) b_2 \quad 
\textrm{and that} \quad
\rho_R( c) = c_1 \left(\id^{\ot n-2} \ot R\right) c_2,
$$
where $ b_i, c_i$ are elements of $A_{R, n-1}\subset A_{R,n}$.
Since $ b_2 c_1$ belongs to $A_{R, n-1} = C_{n-1}$ (by the induction 
hypothesis), we have that
$$
b_2 c_1 =d + \sum_i e_i \left(\id^{\ot n-3} \ot R \ot \id \right) f_i
$$
where $ d, e_i, f_i$ belong to $A_{R, n-2}\subset A_{R,n-1}\subset 
A_{R,n}$. We deduce that
\bea
\rho_R( bc) &=& b_1 \left(\id^{\ot n-2} \ot R\right) d \left(\id^{\ot 
n-2} \ot R\right) c_2\\
&&+\sum _i b_1 \left(\id^{\ot n-2} \ot R\right) e_i
\left(\id^{\ot n-3} \ot R \ot \id \right) f_i  \left(\id^{\ot n-2} 
\ot R\right) c_2\\
&=& b_1 d \left(\id^{\ot n-2} \ot R^2\right) c_2\\
&&+\sum_i b_1 e_i \left(\id^{\ot n-2} \ot R\right) \left(\id^{\ot 
n-3} \ot R \ot \id \right)
\left(\id^{\ot n-2} \ot R\right)  f_i c_2.
\eea
Since the minimal polynomial of $R$ is quadratic,
we have that
$$b_1 d \left(\id^{\ot n-2} \ot R^2\right) c_2 \in C_n.$$
Also, it follows from the braid equation that
\bea
&&b_1 e_i \left(\id^{\ot n-2} \ot R\right) \left(\id^{\ot n-3} \ot R 
\ot \id \right)
\left(\id^{\ot n-2} \ot R\right)  f_i c_2\\
&=& b_1 e_i \left(\id^{\ot n-3} \ot R \ot \id\right) \left(\id^{\ot 
n-2} \ot R  \right)
\left(\id^{\ot n-3} \ot R \ot \id\right)  f_i c_2
\in C_n.
\eea
We conclude that $\rho_R( bc)$ indeed belongs to $C_n$.
\end{proof}

\begin{proof}[Proof of Proposition \ref{prop:quadratic}]
By the previous lemma, any $ x \in A_{R,n}$ can be written as
$$
x = y + \sum_i z_i \left(\id^{\ot n-2} \ot R\right) t_i
$$
with $y, z_i, t_i \in A_{R,n-1} \subset A_{R,n}$. So, we obtain that
\bea
\hbox{Sp}_n \left(x \left(\id^{\ot n-1} \ot \mu \right)\right) &=&
\hbox{Sp}_n \left(y \left(\id^{\ot n-1} \ot \mu \right)\right)+\\
&&
\sum_i \hbox{Sp}_n \left(z_i \left(\id^ {\ot n-2} \ot R\right) t_i 
\left(\id^{\ot n-1} \ot \mu \right)\right)\\
&=&
\hbox{Sp}_n \left(y \ot \mu \right) + \sum_i
\hbox{Sp}_n \left( z_i \left(\id^{\ot n-2} \ot \left(R \circ (\id \ot 
\mu)\right)\right)t_i\right)\\
&=&
\hbox{Sp}_n \left(y \ot \mu \right) + \sum_i
z_i \hbox{Sp}_n \left(\id^{\ot n-2} \ot \left(R \circ (\id \ot 
\mu)\right)\right)t_i\\
&=&
y \ot \hbox{Sp}_1\left(\mu \right) + \sum_i
z_i\left(\id^{\ot n-2} \ot \hbox{Sp}_2\left(R \circ (\id \ot 
\mu)\right)\right)t_i\\
&=& \hbox{Trace}( \mu )\cdot y + \alpha \beta\cdot \sum_i z_i t_i.
\eea
In these identities, we have used elementary properties of the 
operator trace and the condition
of enhancement $\left(E_2^+\right)$.
We conclude that
$\hbox{Sp}_n \left(x \left(\id^{\ot n-1} \ot \mu \right)\right)$ 
belongs to $A_{R,n-1}$.
Thus, $ S $ is redundant.
\end{proof}

\section{Invariants of oriented links derived from algebra structures}

\label{sec:invariants_associated_to_algebras}

In this section, we consider a (unitary associative) $k$--algebra $A$ 
with finite dimension at least $2$,
as well as nonzero scalars $x$ and $y$.
Let $R_{x,y}: A\ot A \to A\ot A$ be the YB operator from \S 
\ref{sec:YB} defined by
$$
R_{x,y}(a\ot b)=x \cdot ab \ot 1 + y \cdot 1 \ot ab - x \cdot a \ot b.
$$
We apply Turaev's procedure to this Yang--Baxter operator.

\subsection{Enhancements of $R_{x,y}$}

Let $ \mu :A \rightarrow A $ be a bijective $k$--linear map
and let $\alpha,\beta$ be nonzero scalars. Let us look for necessary 
and sufficient conditions
for $\left( R_{x,y} , \mu , \alpha , \beta \right)$ to be an enhanced 
YB operator
in the sense of Definition \ref{def:enhancement}.\\

Firstly, condition $\left(E_1\right)$ means that $\mu$ is an automorphism
of the YB operator $R_{x,y}$ (in the sense of Definition \ref{def:isoYB}).
Thus, by Remark \ref{rem:isoYB}, $\mu$ has this property if and only if
there exists a  scalar $c$ such that $c\cdot \mu$ is a $k$--algebra 
automorphism.
Since $\left(R_{x,y},\mu,\alpha,\beta\right)$ is an enhanced
$YB$ operator if and only if $\left(R_{x,y},c \cdot \mu,\alpha,c 
\beta\right)$ is,
we can assume presently with no loss of generality that $\mu$ is an 
automorphism of the algebra $A$.\\

Secondly, let us find an equivalent statement for condition 
$\left(E_2^+\right)$.
For this, we fix a basis $e=( e_1 , e_2 , \dots , e_d )$ of $A$
such that $e_1 = 1$, the unit element of $A$.
Setting $f:= R_{x,y} \circ ( \id \ot \mu )$, we wish to compute 
Sp$_2(f)$ which is defined by
\begin{equation}
\label{eq:Sp2}
\hbox{Sp}_2(f)(e_i)=\sum_{1\leq j,l\leq d}f^{jl}_{il}e_j.
\end{equation}
In the sequel, we will denote by $\rho_{il}\in k$ the $l$--th coordinate
of $e_i \mu(e_l)$ in the basis $e$.

Since $\mu(1)=1$, we have that $f(1 \ot 1)= y\cdot 1 \ot 1$, and so
\begin{equation}  \label{a}
\left\{\begin{array}{ll}
f^{11}_{11}=y, &\\
f^{j1}_{11}=0 & \textrm{ if } j>1.
\end{array}\right.
\end{equation}
Since $ f(1 \ot e_l)=x\cdot \mu(e_l) \ot 1 + (y-x)\cdot 1 \ot 
\mu(e_l) $, we deduce that
\begin{equation}  \label{b}
\left\{\begin{array}{ll}
f^{1l}_{1l}=(y-x) \mu^l_l &\textrm{ if } l>1,\\
f_{1l}^{jl}=0 & \textrm{ if } l>1\textrm{ and } j>1.
\end{array}\right.
\end{equation}
Since $f(e_i \ot 1)= y\cdot 1 \ot e_i$, we deduce that
\begin{equation}  \label{c}
f^{j1}_{i1}=0 \ \ \textrm{ if } i>1 \textrm{ and for any } j.
\end{equation}
Finally, since $f(e_i \ot e_l)= x\cdot e_i \mu(e_l) \ot 1 + y\cdot 1 
\ot e_i \mu(e_l) - x\cdot e_i \ot \mu( e_l)$, we deduce that
\begin{equation}  \label{d}
\left\{\begin{array}{ll}
f^{1l}_{il}= y \rho_{il}& \textrm{ if } i>1 \textrm{ and } l>1,\\
\ f^{il}_{il}=-x \mu^l_l& \textrm{ if } i>1 \textrm{ and } l>1,\\
f^{jl}_{il}=0& \textrm{ if } i>1, l>1 \textrm{ and } j\neq 1,i.
\end{array}\right.
\end{equation}

 From equations (\ref{a},\ref{b},\ref{c},\ref{d}) and formula 
(\ref{eq:Sp2}), we deduce that
$$
\left\{\begin{array}{ll}
\hbox{Sp}_2(f)(e_1)= \left(x + (y-x)\hbox{Trace}( \mu )\right)\cdot e_1,&\\

\hbox{Sp}_2(f)(e_i)= y\left(\sum_{l>1} \rho_{il}\right)\cdot e_1-
x\left(\hbox{Trace}( \mu )-1\right)\cdot e_i & \textrm{ if } i>1.
\end{array}\right.
$$
Hence, we conclude that
\bea
\left(E_2^+\right) & \Longleftrightarrow &
\left\{\begin{array}{ll}
x+(y-x)\hbox{Trace}(\mu) = \alpha \beta & \\
y \sum_{l>1} \rho_{il} = 0& \textrm{ for all } i>1\\
-x\left(\hbox{Trace}(\mu)-1\right)= \alpha \beta &
\end{array}\right.\\
&\Longleftrightarrow&\left\{\begin{array}{ll}
x=\alpha \beta&\\
\hbox{Trace}(\mu)=0&\\
\sum_{l>1} \rho_{il} = 0& \textrm{ for all } i>1
\end{array}\right.\\
&\Longleftrightarrow&\left\{\begin{array}{ll}
x=\alpha \beta&\\
\hbox{Trace}\left(A\to A, b \mapsto a \mu(b)\right)=0 & \textrm{ for 
all } a\in A.
\end{array}\right.
\eea
In the last equivalence, we use the fact that, for any $i>1$, the sum
$\sum_{l>1} \rho_{il}$ coincides with the trace of the $k$--endomorphism
of $A$ defined by $b \mapsto e_i \mu(b)$.\\

As for condition $\left(E_2^-\right)$, we can proceed similarly
to compute Sp$_2\left({R_{x,y}}^{-1}\circ (\id\otimes \mu)\right)$. Recall
from \S \ref{sec:YB} that the inverse of $R_{x,y}$ is given by
$$
{R_{x,y}}^{-1}(a \ot b) = y^{-1}\cdot ab \ot 1 + x^{-1}\cdot 1 \ot ab 
- x^{-1}\cdot a \ot b
$$
so that the same kind of arguments apply. We find that
$$
\left(E_2^-\right) \Longleftrightarrow
\left\{\begin{array}{ll}
y^{-1}=\alpha^{-1} \beta&\\
\hbox{Trace}\left(A\to A, b \mapsto a \mu(b)\right)=0 & \textrm{ for 
all } a\in A.
\end{array}\right.
$$
\vspace{0.1mm}

The above discussion can be summed up as follows.
\begin{lm}
\label{lem:enhancement}
Let $\mu: A \to A$ be a bijective $k$--linear map and let 
$\alpha,\beta$ be nonzero scalars.
Then, $\left(R_{x,y},\mu,\alpha,\beta\right)$ is an enhanced YB 
operator if, and only if,
there exists a scalar $c$ such that $c \cdot \mu$ is a $k$--algebra 
automorphism
and the following conditions are satisfied
\begin{equation}
\label{eq:enhancement}
\left\{\begin{array}{ll}
x = c \alpha \beta &\\
y^{-1}=c \alpha^{-1} \beta&\\
{\rm Trace}\left(A\to A, b \mapsto a \mu(b)\right)=0 & \hbox{ for all } a\in A.
\end{array}\right.
\end{equation}
\end{lm}

In fact, it is enough to consider a
two-dimensional algebra to find a solution to (\ref{eq:enhancement}).

\begin{ex}
\label{ex:2-dimensional}
Let us consider the field
$$
k=\Q\left(x^{1/2},y^{1/2}\right)
$$
of rational fractions in $x^{1/2}$ and $y^{1/2}$, and the $k$--algebra
$$
A =\frac{k[t]}{\left(t^2\right)}.
$$
Then, the automorphism $\mu$ of the algebra $A$ defined by 
$\mu(t)=-t$, the scalars $\alpha=x^{1/2}y^{1/2}$
and $\beta=x^{1/2}y^{-1/2}$ satisfy conditions (\ref{eq:enhancement}) 
with $c=1$, so that
$\left(R_{x,y},\mu,\alpha,\beta\right)$ is an enhanced YB operator.
\end{ex}

\subsection{Connection with the Alexander polynomial}

Assume that $S=\left(R_{x,y},\mu,\alpha,\beta\right)$ is an enhanced 
YB operator.
Lemma \ref{lem:enhancement} gives us necessary and sufficient 
conditions for that:
in particular, by the third condition of (\ref{eq:enhancement}), we 
must have Trace$(\mu)=0$.
That enhanced YB operator is redundant since the minimal polynomial
of $R_{x,y}$ has degree $2$  (Proposition \ref{prop:quadratic}).
Thus, by Remark \ref{rem:trivialize}, the invariant $X_S$ associated 
to $S$ vanishes.
Consequently, we will consider the invariant $X_{S,1}$ associated
to $S$ by taking modified Markov traces of braids.

\begin{te}
\label{th:conclusion} Let $\mu: A \to A$ and $\alpha,\beta \in k$ be such that
$S=\left(R_{x,y},\mu,\alpha,\beta\right)$ is an enhanced YB operator.
Then, the invariant of oriented links $X_{S,1}$ is determined by the 
Alexander polynomial.
Conversely, the Alexander polynomial can be recovered from $X_{S,1}$ 
for appropriate ground fields $k$.
\end{te}

Before proving Theorem \ref{th:conclusion}, we recall that the 
\emph{Alexander polynomial}
of an oriented link $L$, denoted by
$$
\Delta(L)\in \Z\left[t^{\pm 1/2}\right],
$$
is a classical invariant. Defined from the homology of the maximal abelian
covering of the complement of the link (see, for instance, \cite{R}),
the Alexander polynomial also satisfies the following \emph{skein} relations:
\begin{equation}
\label{eq:skein}
\left\{\begin{array}{l}
\Delta\left(\unknot\right)= 1\\
  \Delta \left(\overcrossing\right) - \Delta \left(\undercrossing\right) =
\left(t^{1/2}-t^{-1/2}\right) \Delta \left(\smoothing\right).
\end{array}\right.
\end{equation}
In fact, these relations determine the Alexander polynomial since any link
can be transformed to the trivial link (with the same number of components)
by changing finitely many crossings.

\begin{proof}[Proof of Theorem \ref{th:conclusion}]

Suppose that $\mu:A\to A$ and $\alpha,\beta\in k$ are such that
$S=\left(R_{x,y},\mu,\alpha,\beta\right)$ is an enhanced YB operator
(as in Example \ref{ex:2-dimensional}, for instance).
Then, by Lemma \ref{lem:enhancement}, there exists a nonzero scalar $c$
such that conditions (\ref{eq:enhancement}) holds. In particular, we 
have that $\alpha^2=xy$.
 From (\ref{eq:quadratic}), we deduce that
$$
\alpha^{-1}\cdot R_{x,y} - \alpha\cdot {R_{x,y}}^{-1}=
\left(c^{-1}\beta^{-1}-c\beta\right)\cdot \id^{\ot 2}.
$$
By definition of $X_{S,1}$ and, in particular, by linearity of the maps
Sp$_n$'s, we obtain that
$$
X_{S,1}\left(\overcrossing\right) - X_{S,1}\left(\undercrossing\right) =
\left(c^{-1}\beta^{-1}-c\beta\right)\cdot X_{S,1}\left(\smoothing\right).$$
Moreover, we have that
$$
X_{S,1}\left(\unknot\right) = \beta^{-1}.
$$
Thus, we deduce from (\ref{eq:skein}) that, for any oriented link $L$,
$X_{S,1}$ can be computed from the Alexander polynomial as follows:
\begin{equation}
\label{eq:relation}
X_{S,1}(L)=\beta^{-1} \cdot \Delta(L)\vert_{t^{1/2}=c^{-1}\beta^{-1}}.
\end{equation}

Conversely, let us add the hypothesis that
the ground field $k$ is such that the ring homomorphism
$\Z\left[t^{\pm 1/2}\right] \to k$ defined by $t^{1/2}\mapsto 
c^{-1}\beta^{-1}$ is injective.
(For instance, $k=\Q\left(x^{1/2},y^{1/2}\right)$ in Example 
\ref{ex:2-dimensional} works.)
Then, relation (\ref{eq:relation}) determines $\Delta(L)$ from $X_{S,1}(L)$.

\end{proof}

\textbf{Acknowledgments.}
The authors wish to thank S. Papadima for discussions on that subject 
and comments on that paper.
The second author has been supported by a Marie Curie Fellowship 
(HPMF-CT-2002-01782)
at the University of Wales, Swansea. The first author has been 
supported by a EURROMMAT Fellowship
at the Institute of Mathematics of the Romanian Academy.

\bibliographystyle{amsalpha}

\begin{thebibliography}{99}
\newcommand{\no}{$\textrm{n}^{\circ}$}

\bibitem{DN}  S. D\u asc\u alescu, F. F. Nichita, \emph{Yang--Baxter operators
arising from (co)algebra structures}, Comm. Algebra \textbf{27} \no 12 (1999),
5833--5845.

\bibitem{HOMFLY} P. Freyd, D. Yetter, J. Hoste, W. B. R. Lickorish, 
K. Millett, A. Ocneanu,
\emph{A new polynomial invariant of knots and links},
Bull. Amer. Math. Soc. \textbf{12} (1985), 239--246.

\bibitem{Jimbo} M. Jimbo, \emph{Quantum R matrix for the generalized Toda
system}, Lett. Math. Phys. \textbf{11} (1986), 247--252.

\bibitem{Jones} V. F. R. Jones, \emph{A polynomial invariant for 
knots via von Neumann algebras},
Bull. Amer. Math. Soc. \textbf{12} (1985), 103--111.

\bibitem{Kassel} C. Kassel, Quantum Groups,
Graduate Texts in Math. \textbf{155}, Springer Verlag, 1995.

\bibitem{K} L. H. Kauffman, \emph{New invariants in the theory of knots},
Ast\'erisque \textbf{163}-\textbf{164} (1988), 137--219.

\bibitem{LR} L. Lambe, D. Radford,  Introduction to the Quantum
Yang--Baxter Equation and Quantum Groups: an Algebraic Approach, Math. and its
Applications \textbf{423}, Kluwer Academic Publishers, 1997.

\bibitem{M} J. Murakami, \emph{A state model for the multi--variable 
Alexander polynomial},
Pacific J. Math. \textbf{157} \no 1 (1993), 109--135.

\bibitem{N} F. F. Nichita, \emph{Self-inverse Yang--Baxter operators from
(co)algebra structures}, J. Algebra \textbf{218} \no 2 (1999), 738--759.

\bibitem{PT} J. H. Przytycki, P. Traczyk, \emph{Invariants of links 
of Conway type},
Kobe J. Math. \textbf{4} \no 2 (1987), 115--139.

\bibitem{R} D. Rolfsen, Knots and links, Math. Lecture Series \textbf{7},
Publish or Perish, 1990.

\bibitem{T} V. Turaev, \emph{The Yang--Baxter equation and invariants 
of links},
Invent. Math. \textbf{92} (1988), 527--553.

\end{thebibliography}

\vspace{1.4cm}

\textsc{ \footnotesize G. Massuyeau, Institute of Mathematics of the 
Romanian Academy,
P.O. Box 1--764, 014700 Bucharest, Romania.}

\emph{\footnotesize E-mail address:} \texttt{\footnotesize 
Gwenael.Massuyeau@imar.ro}\\

\textsc{\footnotesize F. F. Nichita, Institute of Mathematics of the 
Romanian Academy,
P.O. Box 1--764, 014700 Bucharest, Romania.}\\
\emph{\footnotesize and} \textsc{\footnotesize Department of 
Mathematics, University of Wales Swansea, Singleton Park,
Swansea SA2 8PP, United Kingdom.}

\emph{\footnotesize E-mail addresses:}
\texttt{\footnotesize Florin.Nichita@imar.ro} \emph{\footnotesize \ and \ }
\texttt{\footnotesize F.F.Nichita@swansea.ac.uk}

\end{document}